\documentclass[10pt]{article}

\usepackage[english]{babel}
\usepackage{ae}
\usepackage{tikz}
\usetikzlibrary{arrows}
\usetikzlibrary{matrix}
\usepackage{amsmath, amsfonts, amssymb, amsthm, color, bm} %mathtools,
\usepackage{url}
\usepackage{authblk}
\usepackage{nicefrac}
\usepackage{fancyhdr}
\usepackage{upgreek}
\usepackage{gensymb}

\usepackage{graphicx,placeins}
\usepackage{subfig}

\newcounter{iabc}

\makeatletter

\newcommand{\rmi}{{\mathrm{i}}}
\newcommand{\rme}{{\mathrm{e}}}

\newcommand{\sign}{{\operatorname{sign\,}}}

\newtheorem{remark}{Remark}

\providecommand{\tabularnewline}{\\}

 \let\oldforeign@language\foreign@language
 \DeclareRobustCommand{\foreign@language}[1]{%
   \lowercase{\oldforeign@language{#1}}}

\makeatother

\begin{document}

\title{Computational Realization of a Non-Equidistant Grid Sampling in Photoacoustics with a Non-Uniform FFT}

%\author{Julian~Schmid,~Thomas~Glatz,~Behrooz~Zabihian,~Mengyang~Liu,~Wolfgang~Drexler~and~Otmar~Scherzer}

\author[1]{Julian Schmid}
\author[1]{Thomas Glatz}
\author[2]{Behrooz~Zabihian}
\author[2]{Mengyang~Liu}
\author[2]{Wolfgang~Drexler}
\author[1,3]{Otmar Scherzer}

\affil[1]{\footnotesize Computational Science Center, University of Vienna, Oskar-Morgenstern-Platz\ 1, 1090 Vienna, Austria}
\affil[2]{\footnotesize Center for Medical Physics and Biomedical Engineering, Medical University of Vienna, Waehringer Guertel\ 18-20, 1090 Vienna, Austria}
\affil[3]{\footnotesize Radon Institute of Computational and Applied Mathematics, Austrian Academy of Sciences, Altenberger Str.\ 69, 4040 Linz, Austria}

\maketitle

\begin{abstract}
\noindent
To obtain the initial pressure from the collected data on a planar sensor arrangement in Photoacoustic tomography, there exists an exact analytic frequency domain reconstruction formula. An efficient realization of this formula needs to cope with the evaluation of the data’s Fourier transform on a non-equispaced mesh. In this paper, we use the non-uniform fast Fourier transform to handle this issue and show its feasibility in 3D experiments. This is done in comparison to the standard approach that uses polynomial interpolation. Moreover, we investigate the effect and the utility of flexible sensor location on the quality of photoacoustic image reconstruction. The computational realization is accomplished by the use of a multi-dimensional non-uniform fast Fourier algorithm, where non-uniform data sampling is performed both in frequency and spatial domain. We show that with appropriate sampling the imaging quality can be significantly improved. Reconstructions with synthetic and real data show the superiority of this method.
\end{abstract}

\section{Introduction}

Photoacoustic tomography is an emerging imaging
technique that combines the good contrast of optical absorption with
the resolution of ultrasound images (see for instance \cite{Wan09}). 
In experiments an object is irradiated by a short-pulsed laser beam. 
Depending on the absorption properties of the material, some light 
energy is absorbed and converted into heat. 
This leads to a thermoelastic expansion, which causes a pressure rise, 
resulting in an ultrasonic wave, called photoacoustic signal. 
The signal is detected by an array of ultrasound transducers outside 
the object. Using this signal the pressure distribution at the time of 
the laser excitation is reconstructed, offering a 3D image proportional 
to the amount of absorbed energy at each position. This is the imaging
parameter of Photoacoustics.
Common measurement setups rely on small ultrasound sensors, which are 
arranged \emph{uniformly} along simple geometries, such as planes, 
spheres, or cylinders covering the specimen of interest
(see for instance \cite{XuWan02a,XuWan02b,XuWan03,Wan09,Bea11} below).

For the planar arrangement of point-like detectors there exist several 
approaches for reconstruction, including numerical algorithms based on 
filtered backprojection formulas and time-reversal algorithms (see for 
instance \cite{XuWan02b,KucKun08,XuWan05,XuWan06}). 

The suggested algorithm in the present work realizes a Fourier inversion formula (see 
\eqref{eq:exact} below) using the \emph{non-uniform fast Fourier transform} (NUFFT).
This method has been designed for evaluation of Fourier transforms at non-equispaced  
points in frequency domain, or non-equispaced data points in spatial, respectively 
temporal domain. The prior is called NER-NUFFT (non-equispaced range-non uniform FFT), 
whereas the latter is called NED-NUFFT (non-equispaced data-non uniform FFT).
Both algorithms have been introduced in \cite{Fou03}. Both NUFFT methods haven 
proven to achieve high accuracy and simultaneously reach the computational efficiency
of conventional FFT computations on regular grids \cite{Fou03}.
This work investigates photoacoustic reconstructions from ultrasound signals recorded at 
\emph{non-equispaced} positions on a planar surface. To the best of our knowledge,
this is a novel research question in Photoacoustics, where the use of regular grids 
is the common choice. For the reconstruction we propose a novel combination of NED- and 
NER-NUFFT, which we call NEDNER-NUFFT, based on the following considerations: 
\begin{enumerate}
 \item The discretization of the analytic inversion formulas (see \eqref{eq:exact}) contains evaluation 
 at non-equidistant sample points in frequency domain. 
 \item In addition, and this comes from the motivation of this paper, we consider evaluation at 
 non-uniform sampling points.
\end{enumerate}
The first issue can be solved by a NER-NUFFT implementation: For {\bf 2D} photoacoustic inversion with 
{\bf uniformly} placed sensors on a measurement {\bf line} such an implementation has been considered 
in \cite{HalSchZan09b}. This method was used for biological photoacoustic imaging in \cite{SchZanHolMeyHan11}. 
In both papers the imaging could be realized in 2D because integrating line detectors \cite{BurHofPalHalSch05,PalNusHalBur09} 
have been used for data recording. In this paper, however, the focus is on 3D imaging, because measurements are taken 
with point sensors. Experimentally, we show the applicability and superiority of the (NED)NER-NUFFT reconstruction formula
in three spatial dimensions, compared with standard interpolation FFT reconstruction.
To be precise, in this paper we conduct \emph{three dimensional} imaging implemented using NEDNER-NUFFT, 
with ultrasound detectors aligned \emph{non-uniformly} on a measurement plane.
To easily assess the effect of a given arrangement, also 2D numerical simulations have been conducted, to support the 
argumentation. We quantitatively compare the results with other computational imaging methods: 
As a reference we use the $k$-wave toolbox \cite{TreCox10} with a standard FFT implementation of the inversion algorithm. 
The NEDNER-NUFFT yields an improvement of the lateral and axial resolution (the latter even by a factor two).

The outline of this work is as follows:

In Section \ref{sec:reconst} we outline the basics of the Fourier reconstruction approach by presenting the underlying Photoacoustic model.
We state the Fourier domain reconstruction formula \eqref{eq:exact} in a continuous setting.
Moreover, we figure out two options for its discretization. 
We point out the necessity of a fast and accurate algorithm for computing the occurring discrete Fourier transforms with non-uniform
sampling points. 

In Section \ref{sec:NUFFT} we briefly explain the idea behind the NUFFT.
We state the NER-NUFFT (subsection \ref{subsec:NER_NUFFT}) and NED-NUFFT (subsection \ref{subsec:NED_NUFFT}) 
formulas in the form we need it to realize the reconstruction on a non-equispaced grid.

In Section \ref{sec:comp1_NER_NUFFT} we discuss the 3D experimental setup.
The NER-NUFFT is compared with conventional FFT reconstruction.
A test chart is used to quantify resolution improvements in comparison to the \emph{k-wave} FFT reconstruction with linear interpolation.
In axial direction this improvement was about 170\% while reducing the reconstruction time by roughly 35\%.

In Section \ref{sec:comp2_NED_NUFFT} we then turn to the NEDNER-NUFFT in 2D with simulated data, 
in order to test different sensor arrangements in an easily controllable environment. 
An equiangular arrangement turns out to yield an improvement of over 40\% compared to the best choice of equispaced sensor arrangement.
Furthermore, we use the insights gained from the 2D simulations to develop an equi-steradian sensor arrangement for our 3D measurements. 
We apply our NEDNER-NUFFT approach on these data and quantitatively compare the outcomes with reconstructions from equispaced data obtained by the NER-NUFFT approach. 
Our results show a significant improvement of the already superior NER-NUFFT.

\section{Numerical Realization of a Photoacoustic Inversion Formula}
\label{sec:reconst}

Let $U \subset \mathbb{R}^d$ be an open domain in $\mathbb R^d$, 
and $\Gamma$ a $d-1$ dimensional hyperplane not intersecting $U$.
Mathematically, photoacoustic imaging consists in solving the operator equation
\[
\mathbf{Q}[f]=p|_{\Gamma\times(0,\infty)}\,,
\]
where $f$ is a function with compact support in $U$ and $\mathbf{Q}[f]$ is the trace on $\Gamma\times(0,\infty)$ of the solution of 
the equation

%\commentO{Zu was brauchst $\Omega$. Es wird jetzt mal die naechten beiden spalten nicht erwaehnt. Es geht 
%wohl um den support von $f$}
%%
%%
\begin{equation*}
 \begin{aligned}
  \partial_{tt} p - \Delta p &=0 \text{ in } \mathbb R^d \times (0,\infty)\,,\\
  p(\cdot,0) &=f(\cdot)  \text{ in } \mathbb R^d\,,\\
  \partial_t p(\cdot,0) &=0  \text{ in } \mathbb R^d\;.
 \end{aligned}
\end{equation*}
In other words, the photoacoustic imaging problem consists in identifying the initial source 
$f$ from measurement data $g=p|_{\Gamma\times(0,\infty)}$. 

An explicit inversion formula for $Q$ in terms of the Fourier transforms of $f$ and $g:=\mathbf Q[f]$ has been found in \cite{XuFenWan02}. 
Let $(\bm x,y)\in \mathbb R^{d-1}\times\mathbb R^+$. 
Assume without loss of generality (by choice of proper basis) that $\Gamma$ is the hyperplane described by $y=0$.
Then the reconstruction reads as follows:
\begin{align}
\label{eq:exact}
\mathbf{F}[f]\left(\bm{K}\right)=
\frac{2K_{y}}{\kappa\left(\bm{K}\right)}\mathbf{F}[\mathbf{Q}f]\left(\bm K_{\bm x},\kappa\left(\bm{K}\right)\right).
\end{align}
where $\mathbf{F}$ denotes the $d$-dimensional Fourier transform:
\begin{align*}
\mathbf{F}[f]\left(\bm{K}\right):=\frac1{(2\pi)^{n/2}}\int\limits_{\mathbb{R}^{d}}\rme^{-\rmi\bm{K}\cdot(\bm{x},y)}f(\bm{x})\mathrm{d}\bm{x}\,,
\end{align*}
and
\begin{align*}
\kappa\left(\bm{K}\right)&=\mathrm{sign}\left(K_{y}\right)\sqrt{\bm{K_{x}}^{2}+K_{y}^{2}}\,,\\
\bm{K}&=(\bm {K_{x}},K_{y})\;.
\end{align*}
Here, the variables $\bm x,\bm {K_x}$ are in $\mathbb R^{d-1}$, whereas $y,K_y\in\mathbb R$.

For the numerical realization these three steps have to be realized in discrete form:
We denote evaluations of a function $\varphi$ at sampling points 
$(\bm x_m,y_n)\in (-X/2,X/2)^{d-1}\times(0,Y)$ by 
\begin{equation}\label{eq:eval_general}
\varphi_{m,n}:=\varphi(\bm x_m,y_n)\;.
\end{equation}
For convenience, we will modify this notation in case of evaluations 
on an equispaced Cartesian grid. We define the $d$-dimensional grid
\begin{align*}%\label{eq:equi_grid}
\mathbf G_x\times \mathrm G_y:=\{-N_x/2,\dots,N_x/2-1\}^{d-1}\times \{0,\dots,N_y-1\}\,,
\end{align*}
and assume our sampling points to be located on $\bm m\Delta_x,n \Delta_y$, 
where
\[
(\bm m,n)\in \mathbf G_x\times \mathrm G_y\,,
\]
and write
\begin{equation}\label{eq:eval_equi}
\varphi_{\bm m,n}=\varphi(\bm m\Delta_x,n\Delta_y)\,,
\end{equation}
where $\Delta_x:=X/N_x$ resp. $\Delta_y:=Y/N_y$ are the occurring step sizes.

In frequency domain, we have to sample symmetrically with respect to $K_y$. Therefore, we also introduce the interval
\[\mathrm G_{K_y}:=\{-N_y/2,\dots,N_y/2-1\}.\]
Since we will have to deal with evaluations that are partially in-grid, partially not necessarily in-grid,
we will also use combinations of \eqref{eq:eval_general} and \eqref{eq:eval_equi}.
In this paper, we will make use of discretizations of the source function $f$, 
the data function $g$ and their Fourier transforms $\hat f$ resp. $\hat g$.

Let in the following 
\begin{align*}
\hat{f}_{\bm j,l}\,=\, \sum_{(\bm m,n)\in \mathbf G_x\times \mathrm G_y}f_{\bm m,n}\rme^{-2\pi\rmi(\bm j\cdot\bm m+ln)/(N_x^{d-1}N_y)}
\end{align*}
denote the $d$-dimensional discrete Fourier transform with respect to space and time. 
By discretizing formula \eqref{eq:exact} via Riemann sums it follows
\begin{align}\label{eq:equiv_rec}
\begin{aligned}
\hat{f}_{\bm j,l}\,\approx &\,\frac{2l}{\kappa_{\bm j,l}}\sum\limits_{n\in \mathrm G_y}\rme^{-2\pi \rmi\,\kappa_{\bm j,l}n/N_y}\\
&\cdot\underset{\bm m\,\in\,\mathbf G_x}{\sum} \rme^{-2\pi \rmi(\bm j\cdot\bm m + l n)/N_x^{d-1}}g_{\bm m,n}\,,
\end{aligned}
\end{align}
where
\begin{align*}
\hfill\kappa_{\bm j,l}&=\sign (l)\sqrt{\bm j^2+l^2}\,,\\
(\bm j,l)&\in\mathbf G_x\times\mathrm G_{K_y}\;.
\end{align*}
This is the formula from \cite{HalSchZan09b}.
\begin{remark}
Note that we use the interval notation for the integer multi-indices for notational convenience.
Moreover, we also choose the length of the Fourier transforms to be equal to $N_x$ in the first $d-1$ dimensions, 
respectively.  
This could be generalised without changes in practice.  
\end{remark}
Now, we assume to sample $g$ at $M$, not necessarily uniform, points 
$\bm x_m\in (-X/2,X/2)^{d-1}$:
Then, 
\begin{align}\label{eq:disc_rec}
\begin{aligned}
\hat{f}_{\bm j,l}\,\approx &\, \frac{2l}{\kappa_{\bm j,l}}\underset{n\in\mathrm G_y}{\sum}e^{-2\pi \rmi\kappa_{\bm j,l}n/N_y}\\
&\cdot\underset{m=1}{\overset{M}{\sum}}\; \frac{h_m}{\Delta_x^{d-1}} e^{-2\pi \rmi(\bm j\cdot \bm x_m)/M}g_{m,n}\;.
\end{aligned}
\end{align}
The term $h_m$ represents the area of the detector surface around $\bm x_m$ 
and has to fulfil $\underset{m=1}{\overset{M}{\sum}}h_m=(N_x\Delta_x)^{d-1}=X^{d-1}$.
Note that the original formula \eqref{eq:equiv_rec} can be received from \eqref{eq:disc_rec} by choosing $\{\bm x_m\}$
to contain all points on the grid $\Delta_x \mathbf G_x$.
  
Formula \eqref{eq:disc_rec} can be interpreted as follows: Once we have computed
the Fourier transform of the data and 
evaluated the Fourier transform at non-equidistant points with respect to the third coordinate, 
we obtain the (standard, equispaced) Fourier coefficients of $f$. 
The image can then be obtained by applying standard FFT techniques.

The straightforward evaluation of the sums on the right hand side of \eqref{eq:disc_rec} would lead to a computational 
complexity of order $N_y^2\times M^2$. Usually this is improved by the use of FFT methods, which have the drawback 
that they need both the data and evaluation grid to be equispaced in each coordinate.   
This means that if we want to compute \eqref{eq:disc_rec} efficiently, 
we have to interpolate both in domain- and frequency space. 
A simple way of doing that is by using polynomial interpolation.
It is used for photoacoustic reconstruction purposes for instance in the
\emph{k-wave} toolbox for Matlab \cite{TreCox10}.
Unfortunately, this kind of interpolation seems to be sub-optimal for Fourier-interpolation 
with respect to both accuracy and computational costs \cite{Fou03, XuFenWan02} 

A regularized inverse k-space interpolation has already been shown to yield
better reconstruction results \cite{JaeSchuGerKitFre07}. 
The superiority of applying the NUFFT, compared to linear interpolation,
has been shown theoretically and computationally by \cite{HalSchZan09b}.

\section{The non-uniform fast Fourier transform}
\label{sec:NUFFT}
 
This section is devoted to the brief explanation 
of the theory and the applicability of the non-uniform Fourier transform, 
where we explain both the NER-NUFFT (subsection \ref{subsec:NER_NUFFT}) 
and the NED-NUFFT (subsection \ref{subsec:NED_NUFFT}) in the form
(and spatial dimensions) we utilise them afterwards.

The NEDNER-NUFFT algorithm used for implementing \eqref{eq:disc_rec} essentially (up to scaling factors) consists of the following steps:
\begin{enumerate}
\item
Compute a $d-1$ dimensional NED-NUFFT in the $\bm x$-coordinates due to our detector placement.
\item
Compute a one-dimensional NER-NUFFT in the $K_y$-coordinate as indicated by the reconstruction formula \eqref{eq:disc_rec}.
\item
Compute an equispaced $d$-dim inverse FFT to obtain a $d$ dimensional picture of the initial pressure distribution.
\end{enumerate}
\subsection{The non-equispaced range (NER-NUFFT) case}\label{subsec:NER_NUFFT}
With the NER-NUFFT (non equispaced range -- non-uniform
FFT) it is possible to efficiently evaluate the discrete Fourier transform
at non-equispaced positions in frequency domain. 

To this end, we introduce the one dimensional discrete Fourier transform, evaluated
at non-equispaced grid points $\kappa_l\in\mathbb R$:
\begin{align}\label{eq:NUDFT1}
%\begin{array}{cl}
\hat{\varphi}_{l}=\underset{n\in \mathrm G_y}{\sum}\varphi_{n}\rme^{-2\pi \rmi\kappa_{l}n/N},\quad l=1,\ldots,M.
%\end{array}
\end{align}
In order to find an efficient algorithm for evaluation of \eqref{eq:NUDFT1}, 
we use a window function $\Psi$, an oversampling factor $c>1$ and
a parameter $c<\alpha<\pi(2c-1)$ that satisfy:
\label{en:psi} 
\begin{enumerate}
\item $\Psi$ is continuous inside some finite interval $[-\alpha,\alpha]$
and has its support in this interval and 
\item $\Psi$ is positive in the interval $[-\pi,\pi]$. 
\end{enumerate}
Then (see \cite{Fou03,HalSchZan09b}) we have the following representation 
for the Fourier modes occurring in (\ref{eq:NUDFT1}): 
\begin{align}\label{eq:NUFFT1}
\begin{aligned}
e^{-\rmi x\theta}\,=\,\frac{c}{\sqrt{2\pi}\Psi(\theta)}\sum\limits_{k\in\mathbb{Z}}\hat{\Psi}(x-k/c)\rme^{-\rmi k\theta/c},\; |\theta|\leq\pi\;.
\end{aligned}
\end{align}
By assumption, both $\Psi$ and $\hat \Psi$ are concentrated around $0$. So we approximate the sum over all $k\in\mathbb Z$ 
by the sum over the $2K$ integers $k$ that are closest to $\kappa_l+k$. 
By choosing $\theta=2\pi n/N-\pi$ and inserting \eqref{eq:NUFFT1} in \eqref{eq:NUDFT1}, we obtain
\begin{align}\label{eq:NER_NUFFT}
\begin{aligned}
\hat\varphi_{l}\,&\approx   \,
\sum\limits_{k=-K+1}^K \hat{\Psi}_{l,k}\sum\limits_{n\in \mathrm G_y}
\frac{\varphi_{n}}{\Psi_{n}}e^{-2\pi \rmi ln/cN}\,,\\
l \,& = \, 1, \dots,M\;.
\end{aligned}
\end{align}
Here $K$ denotes the interpolation length and 
\begin{equation*}
\begin{aligned}
 \Psi_n  &:=\Psi(2\pi n/N_y-\pi)\,,\\
 \hat{\Psi}_{l,k}  &:=\frac{c}{\sqrt{2\pi}}\,\rme^{-\rmi\pi(\kappa_l-(\mu_{l,k}))}\hat{\Psi}(\kappa_{l}-(\mu_{l,k}))\,,
\end{aligned}
\end{equation*}
where $\mu_{l,k}$ is the nearest integer (i.e. the nearest
equispaced grid point) to $\kappa_{l}+k$.

The choice of $\Psi$ is made in accordance with the assumptions
above, so we need $\Psi$ to have compact support. Furthermore, to
make the approximation in (\ref{eq:NER_NUFFT}) reasonable, its Fourier
transform $\hat{\Psi}$ needs to be concentrated as much as possible
in $[-K,K]$. In practice, a common choice for $\Psi$ is the Kaiser-Bessel
function, which fulfils the needed conditions, and its Fourier transform
is analytically computable.

\subsection{The non-equispaced data (NED-NUFFT) case}\label{subsec:NED_NUFFT}

A second major aim of the present work is to handle data measured
at non-equispaced acquisition points $\bm x_{m}$ in an efficient and accurate
way. Therefore we introduce the non-equispaced data, $d-1$ dimensional DFT
\begin{align}\label{eq:NUDFT2}
\begin{aligned}
\hat{\varphi}_{\bm j}&=\underset{m=1}{\overset{M}{\sum}}\varphi_{m}\rme^{-2\pi \rmi (\bm j\cdot \bm x_{m})/N}\,,\\
\bm j&\in\mathbf G_x\;.
\end{aligned}
\end{align}
The theory for the NED-NUFFT is largely analogous to the NER-NUFFT \cite{Fou03} 
as described in Subsection \ref{subsec:NER_NUFFT}. 
The representation \eqref{eq:NUFFT1} is here used for each entry of $\bm j$ and inserted 
(with now setting $\theta=2\pi n/N$) into formula \eqref{eq:NUDFT2}, which leads to
\begin{align}\label{eq:NED_NUFFT}
\begin{aligned}%{cl}
\hat\varphi_{\bm j}\,\approx &  
\,\frac1{\Psi_{\bm j}}
\sum\limits_{m=1}^{M}
~\sum\limits_{\bm k\in\{-K,\dots,K-1\}^{d-1}}
\varphi_{m}\hat{\Psi}_{\bm j,\bm k}\\
&\cdot \rme^{-2\pi \rmi \left(\bm j\cdot\bm{\mu}_{m,\bm k}\right)/cM}\,,
\end{aligned}
\end{align}
where the entries in $\bm \mu_{m,\bm k}$ are the nearest integers to $\bm x_m+\bm k$. 
Here we have used the abbreviations
\begin{align*}
\begin{aligned}
\Psi_{\bm j,\bm k} &:=\,\prod\limits_{i=1}^{d-1}\Psi(2\pi \bm j/N_x)\,,\\
\hat\Psi_{\bm j,\bm k} &:=\,\prod\limits_{i=1}^{d-1}\left(\frac{c}{\sqrt{2\pi}}\right)\hat{\Psi}((\bm x_{m})_i-(\bm\mu_{m,\bm k})_i)\,,
\end{aligned}
\end{align*}
for the needed evaluations of $\Psi$ and $\hat{\Psi}$.

Further remarks on the implementation of the NED- and NER-NUFFT, 
as well as a summery about the properties of the Kaiser-Bessel function 
and its Fourier transform can be found in \cite{Fou03,HalSchZan09b}.

\section{Comparison of NER-NUFFT and k-wave FFT}
\label{sec:comp1_NER_NUFFT}

Before we turn to the evaluation of the algorithm we describe the photoacoustic 
setup. Our device consists of a  FP (Fabry P\'{e}rot) polymer film
sensor for interrogation \cite{Bea05,BeaPerMil99}. 
A $50\,\mathrm{Hz}$ pulsed laser source and a subsequent optical parametric 
oscillator (OPO) provide optical pulses. These pulses have a very narrow 
bandwidth and can be tuned within the visible and near infrared spectrum. 
The optical pulses are then transmitted via an optical fibre. 
When the light is emitted it diverges and impinges upon a sample with homogeneous fluence, 
thus generating a photoacoustic signal. 
This signal is then recorded via the FP-sensor head. 
The sensor head consists of an approximately $\mathrm{38\,\upmu m}$ thick polymer
(Parylene C) which is sandwiched between two dichroic dielectric coatings.
These dichroic mirrors have a noteworthy transmission characteristic.
Light from $600$ to $1200\,\mathrm{nm}$ can pass the mirrors largely
unabated, whereas the reflectivity from $1500$ to $1650\,\mathrm{nm}$
(sensor interrogation band) is about 95\% \cite{ZhaLauBea08}. 
The incident photoacoustic wave produces a linear change in the optical
thickness of the polymer film. A focused continuous wave laser, operating
within the interrogation band, can now determine the change of thickness
at the interrogation point via FP-interferometry. 

We choose two 2D targets for comparison, a star and a USAF (US Air Force) resolution test chart. Both
targets are made of glass with a vacuum-deposited durable chromium coating. The star target has 72
sectors on a pattern diameter of 5 mm, with an unresolved core diameter of $100\,\mathrm{\upmu m}$.

The targets are positioned in parallel to the sensor surface
at a distance of about $4\,\mathrm{mm}$, and water is used as coupling
medium between the target and the sensor.

%For all reconstructions an upsampling factor of $2$\commentO{c=2 from the NUFFT, oder ist was anderes gemeint?} is used.
The interpolation length for the
NER-NUFFT reconstruction is $K=6$. The computational
times are shown in Tab. \ref{tab:NER-NUFFT comparison} showing that
the linearly interpolated FFT is about 30 \% slower than the NER-NUFFT. 

\begin{figure} 
\begin{center}\includegraphics[width=0.6\columnwidth]{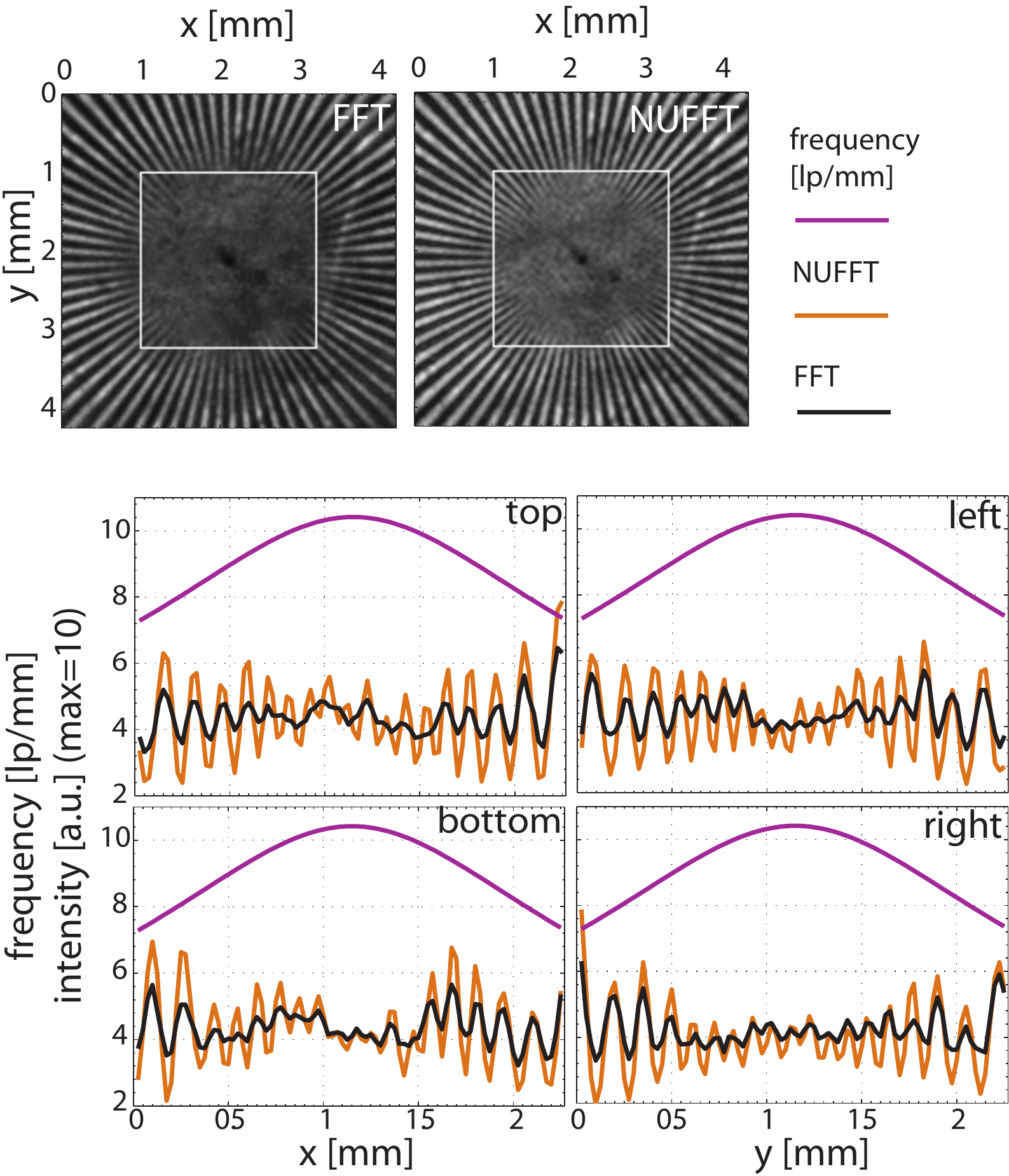}\end{center}
\caption{Segment of the MIP (maximum intensity projection) in the $z$-axis of a star sample, reconstructed
using FFT and linear interpolation (top left) and using the NUFFT (top
right) with a square plotted around the center. The intensity of the
reconstructed image along the sidelines of the square is plotted,
for both reconstructions. The purple line indicates the frequency
of the star sample.\label{fig:StarPhantomReconstruction}}
\end{figure}

A segment of the reconstructed star target is shown in Fig. \ref{fig:StarPhantomReconstruction}.
The intensity is plotted against the sides of an imaginary square ($2.67\,\mathrm{mm}$) placed
around the center of the star phantom. It is clearly visible, that
the FFT reconstruction is not able to represent the line pairs, when
the density exceeds $10\,\mathrm{lp/mm}$, corresponding to a resolution
of $100\,\mathrm{\upmu m}$, whereas they are still largely visible
for the NER-NUFFT reconstruction.

\begin{table*}[tbh]
\caption{Comparison between the NED-NUFFT and FFT reconstruction for a USAF-chart
and a comparison of computational times for both phantoms. The improvement in percent was caluclated by: $100\times(1-\mathrm{FFT/NUFFT})$  \label{tab:NER-NUFFT comparison}}

\begin{center}%
\begin{tabular}{cccc}
\hline 
 & NUFFT & FFT & Improvement\tabularnewline
\hline 
\hline 
FWHM axial LSF & $23.23\pm 0.56\,\mathrm{\upmu m}$ & $62.34\pm 0.62\,\mathrm{\upmu m}$  & $168.47\pm 6.88\,\%$ \tabularnewline
\hline 
FWHM lateral LSF & $33.44\pm 7.95\,\mathrm{\upmu m}$ & $40.82\pm 7.34\,\mathrm{\upmu m}$ & $18.63\pm 8.50\,\%$\tabularnewline
\hline 
Time: Star target & 140 s & 189 s & $35\,\%$\tabularnewline
\hline 
Time: USAF chart & 298 s & 384 s & $29\,\%$\tabularnewline
\hline 
\end{tabular}\end{center}
\end{table*}

For a quantitative comparison of the resolution we use the USAF chart.
It is recorded on an area of $146\times146$ sensor points corresponding
to $1.022\times1.022\,\mathrm{cm^{2}}$ with a grid spacing of $70\,\mathrm{{\upmu m}}$
and a time resolution of $8\,\mathrm{ns}$. 
% A segment of the $xy$-MIP for is shown in Fig. \ref{fig:USAF-Chart}.
As yet there is no standardized procedure to measure the resolution of
a photoacoustic imaging system. We proceed similar to \cite{ZhaLauBea08}, %KimShiRyuSon12
by fitting a \emph{line spread function} (LSF) and an \emph{edge spread function}
(ESF) to the intensities of our reconstructed data. For the LSF to
be meaningful, its source has to approximate a spatial delta function.
This is the case in the $z$-axis, since the chrome coating of the
USAF target is just about $0.1\,\mathrm{\upmu m}$ thick.

\begin{figure}
\begin{center}\includegraphics[width=0.6\columnwidth]{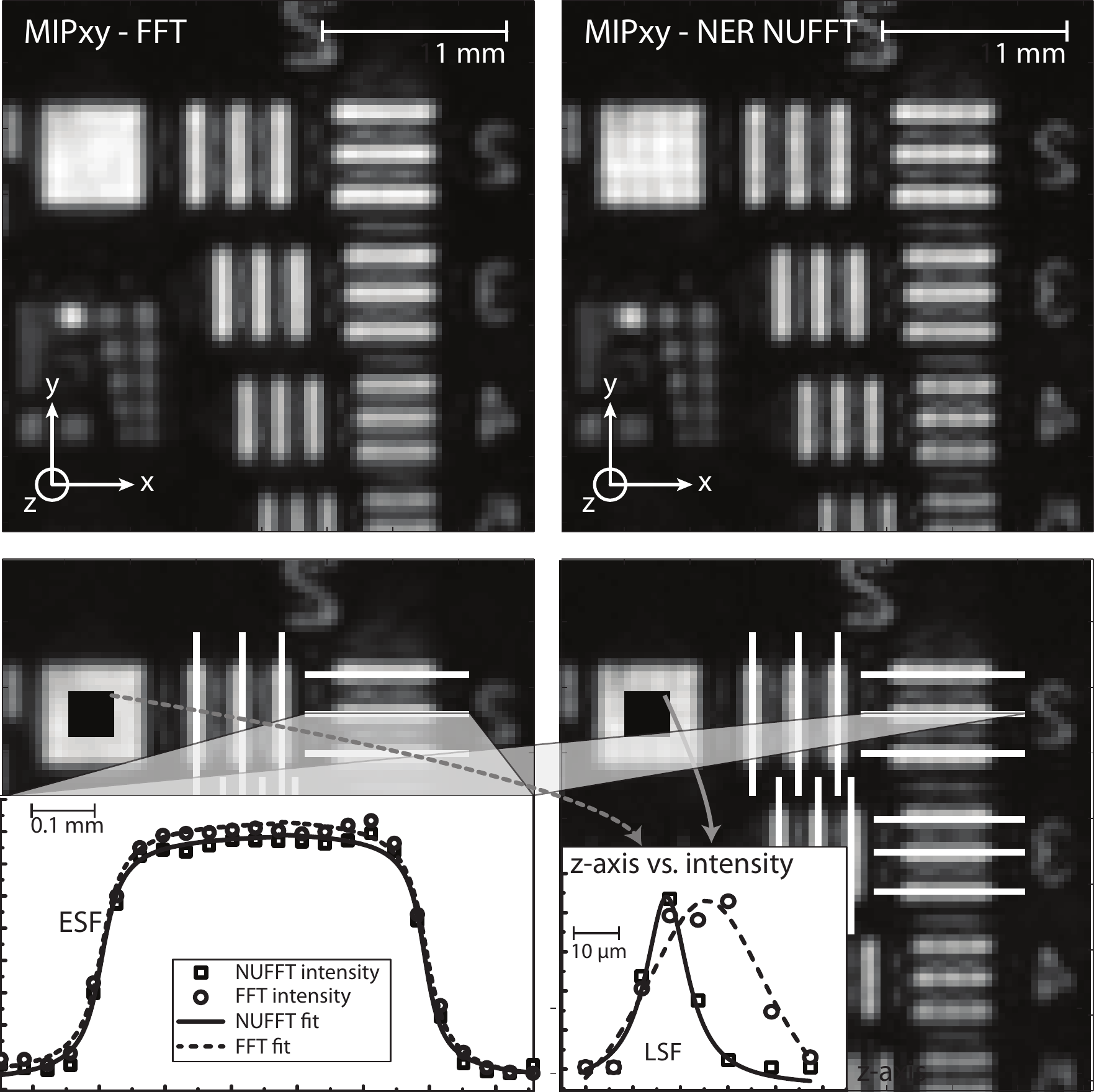}\end{center}

\caption{Segment of the $xy$-MIP of an USAF chart reconstruction conducted
with FFT (left) and NER-NUFFT (right). On the bottom images the data-sets used for the quantitative resolution are marked.
The black square depicts the 49 $xy$-coordinates used for the LSF fit along the $z$-axis for the axial resolution. The data for a single point is shown in the bottom right inlay. 
The white lines show the intensity fit for the lateral resolution. The bottom left inlay shows the data for a single white line.  \label{fig:USAF-Chart}}
\end{figure}

We fit a Cauchy-Lorentz
distribution to the $z$-axis of our reconstructed data for 49 adjacent $xy$-coordinates. Their positions are marked as a black square within the white square, depicted in the bottom images of Fig. \ref{fig:USAF-Chart}. The reconstruction in the $z$-axis for a single point is shown in the bottom right inlay of Fig. \ref{fig:USAF-Chart}, for 8 points, covering a distance of $94.74\,\mathrm{\upmu m}$ in the $z$-direction, around the maximum intensity. A fit of the Lorentz distribution is shown for both reconstruction methods.
The FWHM (\emph{full width, have maximum}) of the Lorentz distribution,
\[I(z)=\frac{2a_{0}w}{\pi\left(w^{2}+4(z-z_{0})^{2}\right)}\,,\]  is the parameter
$w$. The output $I(z)$ is the intensity in dependence from the $z$-axis, and $z_0$ and $a_0$ are fitting parameters.

The average and standard deviation of $w$ for both of the 49 datasets are shown in Table \ref{tab:NER-NUFFT comparison}.
The line spread function FWHM of the FFT-reconstruction turns out to be more than twice as big as the one of the NER-NUFFT reconstruction.

For the lateral resolution, there is no target approximating a delta
function, so we have to use the ESF instead:
\[I(x)=I_{0+}a_{0}\left(\frac{1}{\pi}\arctan\left(\frac{x-x_{0}}{w/2}\right)+\frac{1}{2}\right)\;.\] 
The $w$ here is the FWHM of the associated LSF, and $I_0$, $x_0$ and $a_0$ are fitting parameters. 
The ESF requires a step function as source, of which our target provides plenty. We choose the long side of 12 bars, marked by white lines in the bottom images of Fig. \ref{fig:USAF-Chart}, for this fit. The data for a particular line is shown in the bottom left inlay. We omitted all datasets, where only one point marked the transition from low to high intensity, rendering the edge fit unreliable and resulting in unrealistic improvements of our new method well over 100\%. Finally we averaged over 15 edges. The results are shown in Table \ref{tab:NER-NUFFT comparison}.  While the deviation between different edges is quite high, the improvement of our NER-NUFFT reconstruction for the 15 evaluated edges ranged from $6$ to $44\%$.

\section{Non-equidistant grid sampling}\label{sec:comp2_NED_NUFFT}
%\commentO{Man kann nicht wirklich von sparse grid sprechen, wenn wir immer die gleiche Anzahl von Sensoren verwenden}

%We investigate how non-equidistant sensor positioning effects the imaging quality. 
The current setups allow data acquisition at just one single sensor point for each laser pulse 
excitation. Since our laser is operating at $50\,\mathrm{Hz}$ data recording of a typical 
sample requires several minutes. Reducing this acquisition time is a crucial step in advancing
Photoacoustics towards clinical and preclinical application. 
Therefore, in this work we try to maximize the image quality for a given number of acquisition points. 
We are able to do this, because our newly implemented NEDNER-NUFFT is ideal for dealing with non-equispaced positioned sensors.
This newly gained flexibility of sensor positioning offers many possibilities
to enhance the image quality, compared against a rectangular grid.
For instance a hexagonal grid was found to yield an efficiency of $90.8\,\%$ compared
with $78.5\,\%$ for an exact reconstruction of a wave-number limited
function \cite{PetMid62}.

Also any non-equispaced grids that may arrive from a specific experimental setup can be efficiently computed
via the NEDNER-NUFFT approach.

Here, we will use it to tackle the limited view problem. 
Many papers deal with the limited view problem, when reconstructing images \cite{DanTaoLiuWan12,XuWanAmbKuc04,TaoLiu10}.
Our approach to deal with this problem is different. 
It takes into account that in many cases the limiting factor is the number of sensor
points and the limited view largely a consequence of this limitation.

In our approach we therefore use a grid arrangement that is dense close to a center of interest 
and becomes sparser the further away the sampling points are located. 
We realize this by means of an equiangular, or equi-steradian sensor arrangement,
where for a given point of interest each unit angle or steradian
gets assigned one sensor point.

\begin{figure}
\begin{center}\includegraphics[width=0.6\columnwidth]{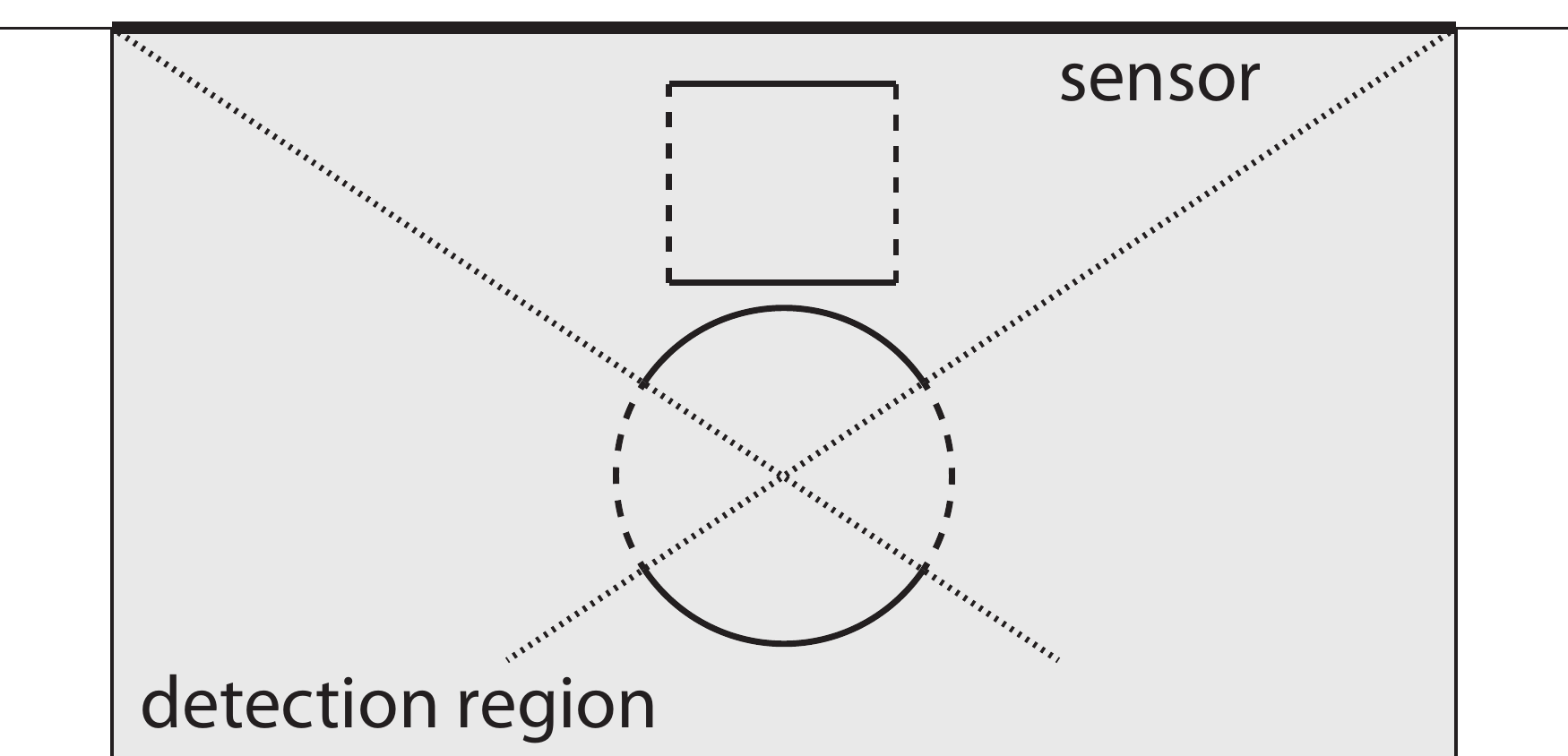}\end{center}

\caption{Depiction of the limited view problem. Edges whose normal vector cannot
intersect with the sensor surface are invisible to the sensor. The
invisible edges are the coarsely dotted lines. The detection region
is marked by a grey background. The finely dotted lines are used to
construct the invisible edges. Edges perpendicular to the sensor surface
are always invisible for a plane sensor. \label{fig:Fieldofview}}
\end{figure}

To understand the limited view problem, it is helpful to define a detection region.
According to \cite{XuWanAmbKuc04}, this is the region which is enclosed
by the normal lines from the edges of the sensor. Pressure waves always
travel in the direction of the normal vector of the boundary of the
expanding object. Therefore certain boundaries are invisible to the
detector as depicted in Fig. \ref{fig:Fieldofview}. 

\subsection{Equiangular and equi-steradian projection sensor mask} 

For the equiangular sensor arrangement a point of interest is chosen.
Each line, connecting a sensor point with the point of interest, encloses
a fixed angle to its adjacent line. In that sense we mimic a circular
sensor array on a straight line. The position of the sensor points
can be seen on top of the third image in Fig. \ref{fig:tree}. 

The obvious expansion of an equiangular projection to 3D is an 
equi-steradian projection. This problem is analogous to the problem of placing
equispaced points on a 3D sphere and then projecting the points, from
the center of the sphere, through the points, onto a 2D plane outside
the sphere. 

We developed an algorithm for this problem, which is explained in detail in Appendix \ref{App:equi-ster}. Our input variables are the grid size, 
the distance of the center of interest from the sensor plane and the desired number of acquisition points, 
which will be rounded to the next convenient value.

A sensor arrangement with 1625 points
on a $226\times226$ grid is shown in the top left image in Fig. \ref{fig:Sensor_masks}.

\subsection{Weighting term} \label{sec:weighting}

To determine the weighting term $h_m$ in Eq. \ref{eq:disc_rec} for 3D we introduce a function
that describes the density of equidistant points per unit area $\rho_{p}$. 
In our specific case, $\rho_{p}$ describes the density on a sphere around a center of interest.
Further we assume that $\rho_{p}$ is spherically symmetric and decreases quadratically with the distance from the center of interest $r$: $\rho_{p,s}\propto1/r^{2}$.
We now define $\rho_{p,m}$ for a plane positioned at distance
$r_{0}$ from the center of interest. In this case $\rho_{p,s}(r)$
attenuates by a factor of $\sin\alpha$, where $\alpha=\arcsin(r_{0}/r)$
is the angle of incidence. Hence $\rho_{p,m}\propto r_{0}/r^{3}$.
%The acquisition point density of the equispaced grid is $\rho_e=1/r_e^2$, where $r_e$ is the grid spacing.
By applying the spacing of the regular grid .
This yields a weighting term of:
%\[h_{m}(r)\propto r^{3}/\left(r_{0}\right)\]
\[h_{m}(r)\propto r^{3}\]

Analogously we can derive $h_m$ for 2D:

\[h_{m}(r)\propto r^{2}\]

We applied a normalization after the reconstruction to all measurements.

For the application of this method to the FP setup it is noteworthy
that there is a frequency dependency on sensitivity which itself depends on the angle of incidence, which has been extensively discussed in \cite{CoxBea07}.
The angle of incidence for our specific setup is $62\degree$.
At this angle, the frequency components around $2\,\mathrm{MHz}$ get attenuated by more than $10\,\mathrm{dB}$.  
Below $1\,\mathrm{MHz}$ the frequency response remains quite stable (attenuation below $5\,\mathrm{dB})$ 
for the measurement angles occurring in our setup.

\section{Computational assessment of different non-equispaced grid arrangements in
two dimensions to tackle the limited view problem}

\begin{figure}
\begin{center}\includegraphics[width=0.6\columnwidth]{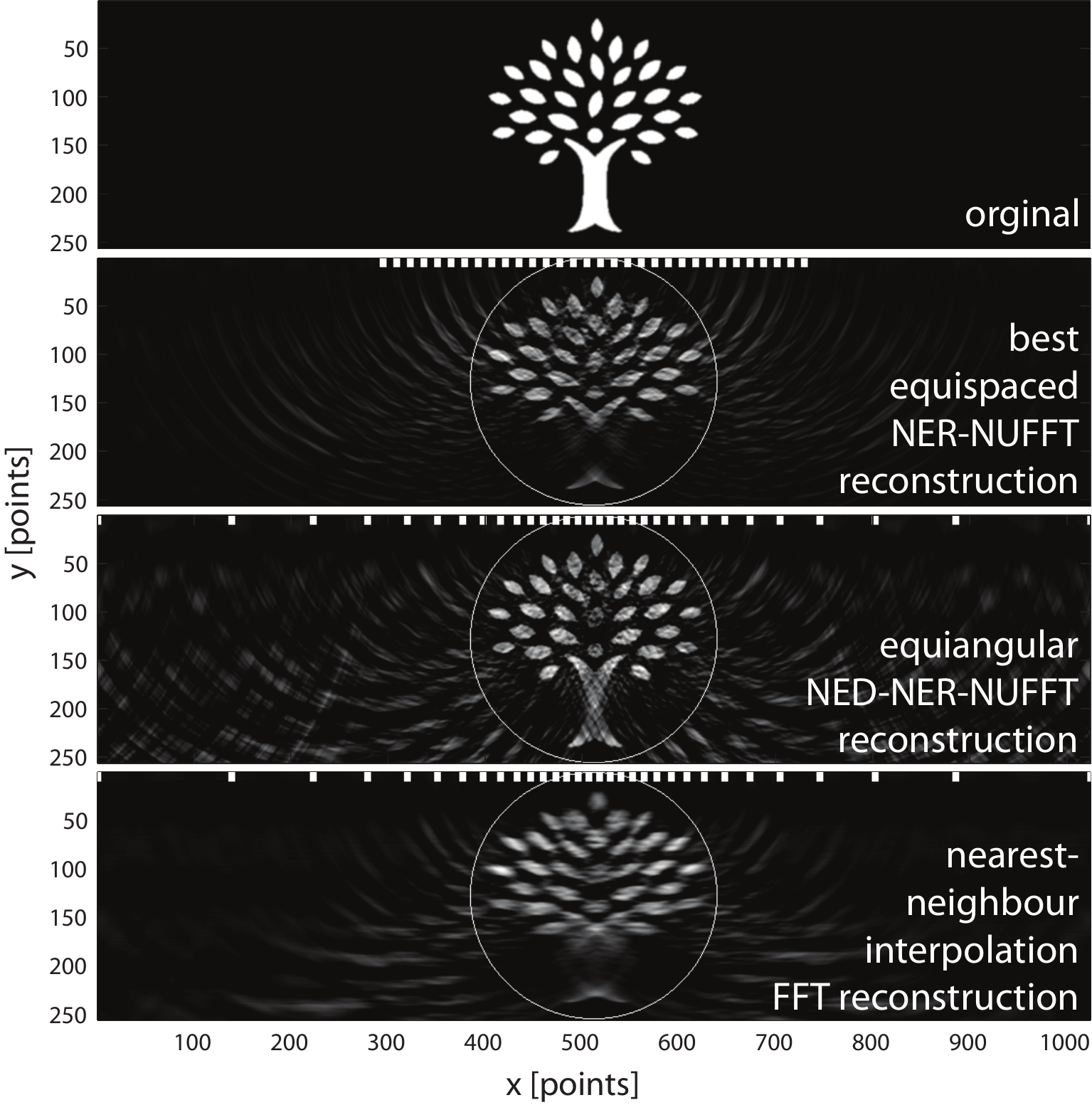}\end{center}\caption{Various reconstructions of a tree phantom (top) with different sensor
arrangements. All sensor arrangements are confined to 32 sensor points.
The sensor positions are indicated as white rectangles on the top
of the images. The second image shows the best (see Fig.\ref{fig:Tree_CC})
equispaced sensor arrangement, with a distance of 13 points between
each sensor. The third image shows the NEDNER-NUFFT reconstruction
with equiangular arranged sensor positions. The bottom image shows
the same sensor arrangement, but all omitted sensor points are linearly
interpolated and afterwards a NER-NUFFT reconstruction was conducted.
\label{fig:tree}}
\end{figure}

A tree phantom, designed by Brian Hurshman and licensed under CC BY
3.0%
\footnote{http://thenounproject.com/term/tree/16622/%
}, is chosen for the 2 dimensional computational experiments on a grid
with $x=1024$ $z=256$ points. A forward simulation is conducted
via \emph{k-wave} \cite{TreCox10}. The forward simulation of the k-wave
toolbox is based on a first order k-space model. A PML (perfectly
matched layer) of 64 gridpoints is added, as well as $30\,\mathrm{dB}$ of 
noise.

In Fig. \ref{fig:tree} our computational phantom is shown at the
top. For each reconstruction a subset of 32 out of the 1024 possible
sensor positions was chosen. In Fig. \ref{fig:tree} their positions
are marked at the top of each reconstructed image. For the equispaced
sensor arrangements, we let the distance between two adjacent sensor
points sweep from 1 to 32. The sensor points where always centered
in the $x$-axis.

To compare the different reconstruction methods we used the correlation coefficient
and the Tenenbaum sharpness. These quality measures are explained in Appendix \ref{App:Qualitmeas}.

\begin{figure}
\begin{center}\includegraphics[width=0.6\columnwidth]{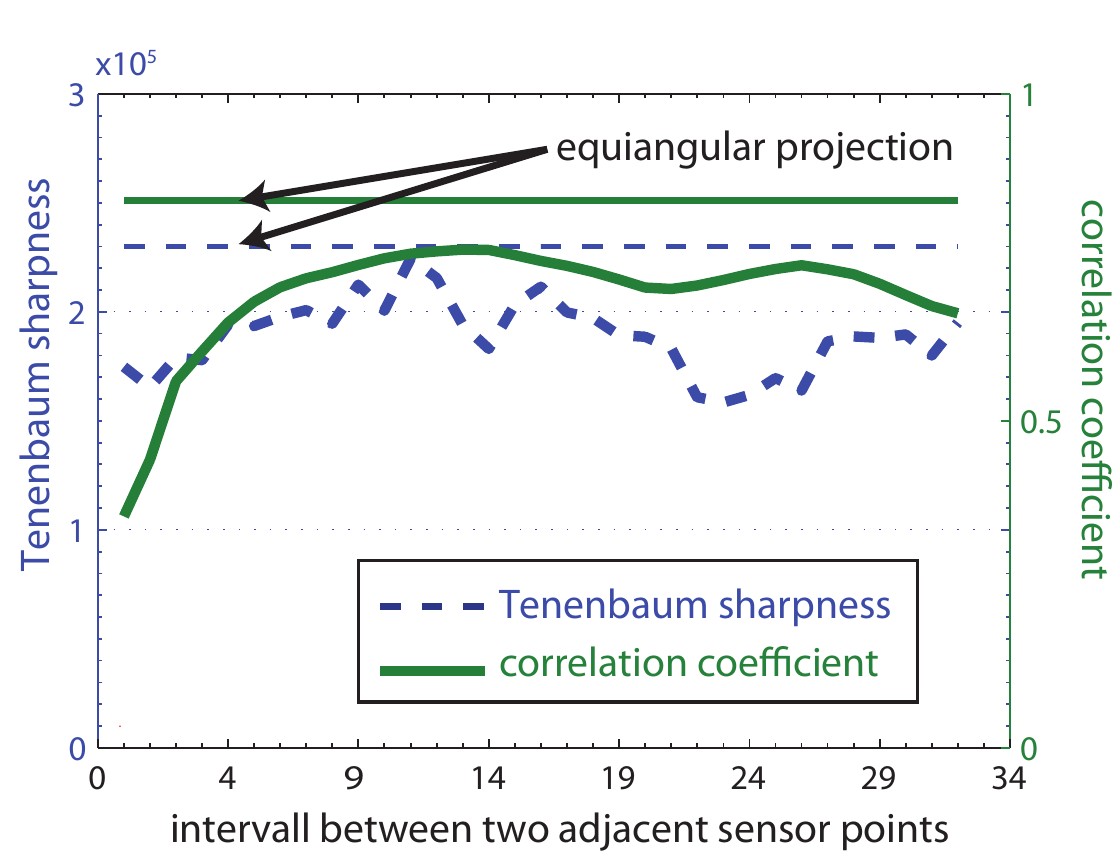}\end{center}
\caption{Correlation coefficient and Tenenbaum sharpness for equispaced sensor
arrangements with intervals between the sensor points reaching from
1 to 32. The maximum of the correlation coefficient is at 13. The
corresponding reconstruction is shown in Fig. \ref{fig:tree}. The
straight lines indicate the results for the equiangular projection.\label{fig:Tree_CC}}
\end{figure}

We applied the correlation coefficient only within
the region of interest marked by the white circle in Fig. \ref{fig:tree}.
The Tenenbaum sharpness was calculated on the smallest rectangle,
containing all pixels within the circle. The results are shown in
Fig. \ref{fig:Tree_CC}. 

The Tenenbaum sharpness for the equiangular sensor placement was
23001, which is above all values for the equispaced arrangements.
The correlation coefficient was 0.913 compared to 0.849, for the best
equispaced arrangement. In other words, the equiangular arrangement
is 42.3 \% closer to a full correlation, than any equispaced grid.

In Fig. \ref{fig:tree} the competing reconstructions are compared.
While the crown of the tree is depicted quite well for the equispaced
reconstruction, the trunk of the tree is barely visible. This is due
to the limited view problem. When the equispaced interval increases,
the tree becomes visible, but at the cost of the crown's quality.
In the equiangular arrangement a trade off between these two effects
is achieved. Additionally the weighting term for the outmost sensors
is 17 times the weighting term for the sensor point closest to the
middle. This amplifies the occurrence of artefacts, particularly
outside of our region of interest. 

The bottom image in Fig. \ref{fig:tree} shows the equiangular
sensor arrangement, reconstructed in a conventional manner. The missing
sensor points are interpolated to an equispaced grid, and a NER-NUFFT
reconstruction is applied afterwards. We conducted a linear interpolation
from our subset to all 1024 sensor points. The correlation coefficient
for this outcome was 0.7348 while the sharpness measure was 21474.
This outcome exemplifies the clear superiority of the NUFFT to conventional
FFT reconstruction when dealing with non-equispaced grids.

\section{3D application of the NED-NER-NUFFT with real data }

\begin{figure}
\begin{center}\includegraphics[width=0.6\columnwidth]{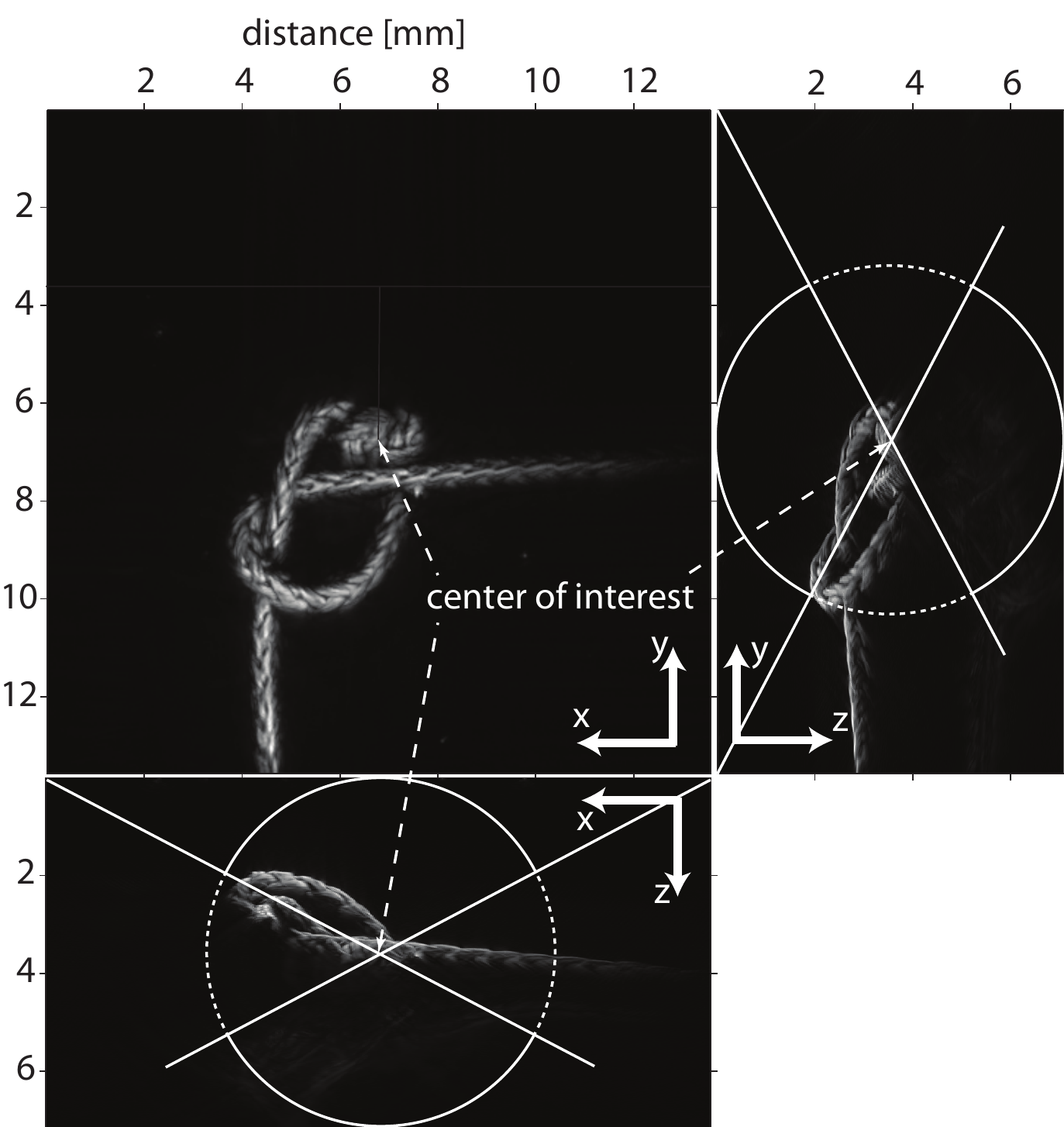}\end{center}

\caption{MIPs of all 3 planes for the NER-NUFFT reconstruction of a yarn phantom.
\label{fig:Golden-Standard}}
\end{figure}

For a qualitative assessment of our new sensor arrangement we need
a 3D phantom. We choose a yarn which we record on a rectangular grid with
$226\times 226=51076$ sensor points, with a grid spacing of $60\,\upmu\mathrm{m}$
and time sampling of $\mathrm{d}t=8\,\mathrm{ns}.$ Hence an area
of $13.56\times 13.56\,\mathrm{mm}$ is covered. As coupling medium water
is used, in which the yarn is fully immersed.

To determine the utility of non-equispaced grid sampling, we follow a certain
routine. First we acquire a densely sampled dataset. Then we use
a very small subset of the initially collected sensor data, to test
different sensor arrangements. Therefore we can always use the complete
reconstruction as our model standard together with the quality measurements explained in Appendix \ref{App:Qualitmeas}.

An upsampling factor of 2 was used for
all reconstructions, hence the reconstructed image of the MIP for
the $xy$ plane consists of $452\times452$ pixels. The complete reconstruction
with the NER-NUFFT took 154 seconds.
Maximum intensity projections (MIP) of this full reconstruction for
all axis are shown in Fig. \ref{fig:Golden-Standard}. 
The MIP in
the $xy$ plane is our model standard, for comparison with the other
reconstructions.

\begin{figure}
\begin{center}\includegraphics[width=0.6\columnwidth]{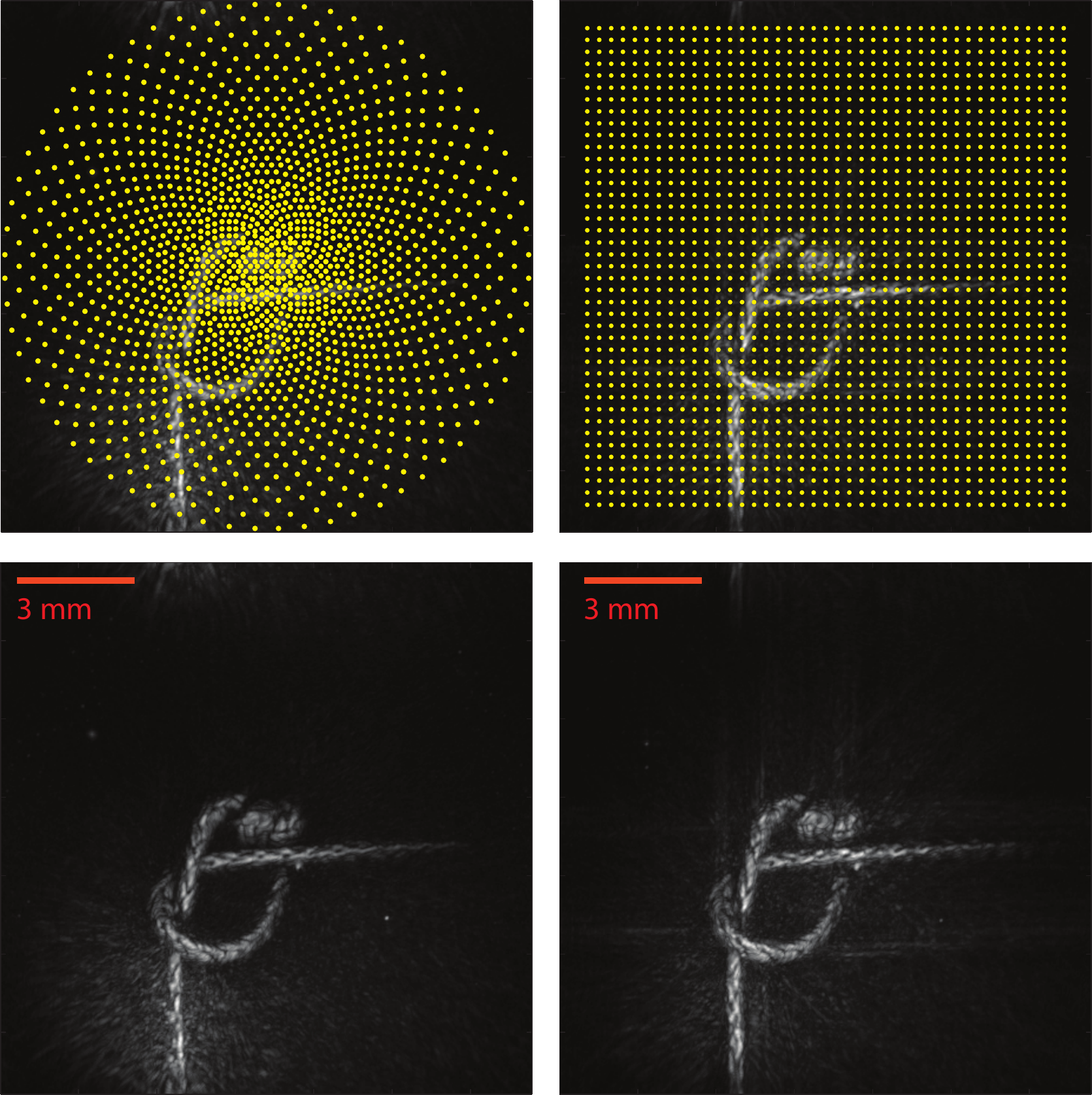}\end{center}

\caption{The sensor placement for the cirular arrangement is shown on the top left
image, comprising 1625 sensor points. On the top right an equispaced sensor arrangement
with $41\times41=1689$ sensor points is displayed. The intervall between
two adjacent sensor points is $5$ for this configuration. The bottom images show
MIPs of the NED-NER-NUFFT reconstruction only using the sensor points shown above.
\label{fig:Sensor_masks}}
\end{figure}

For the equi-steradian sensor mask we choose our center of interest right in
the center of the $xy$-MIP where the little knot can be seen, $3.6\,\mathrm{mm}$
off the sensor surface. The resulting sensor mask, including the reconstruction
is shown in Fig. \ref{fig:Sensor_masks}, it consists of $1625$ sensor
points (or 3.18 \% of the initial number of sensor points). The weighting
term, accounting for sensor sparsity is 9.7 times higher for the outermost
sensor point, than for centermost sensor points in the $xy$-plane. The reconstruction
for this arrangement took 134 seconds.

We compared this arrangement to rectangular grids, which all had $41\times41=1681$ sensor points (or 3.29 \% of the initial
number of sensor points) but varying distances between two adjacent
points. The grid with an interval between sensor points of
5 is shown on the top left in Fig. \ref{fig:Sensor_masks}. 

\begin{figure}
\begin{center}\includegraphics[width=0.6\columnwidth]{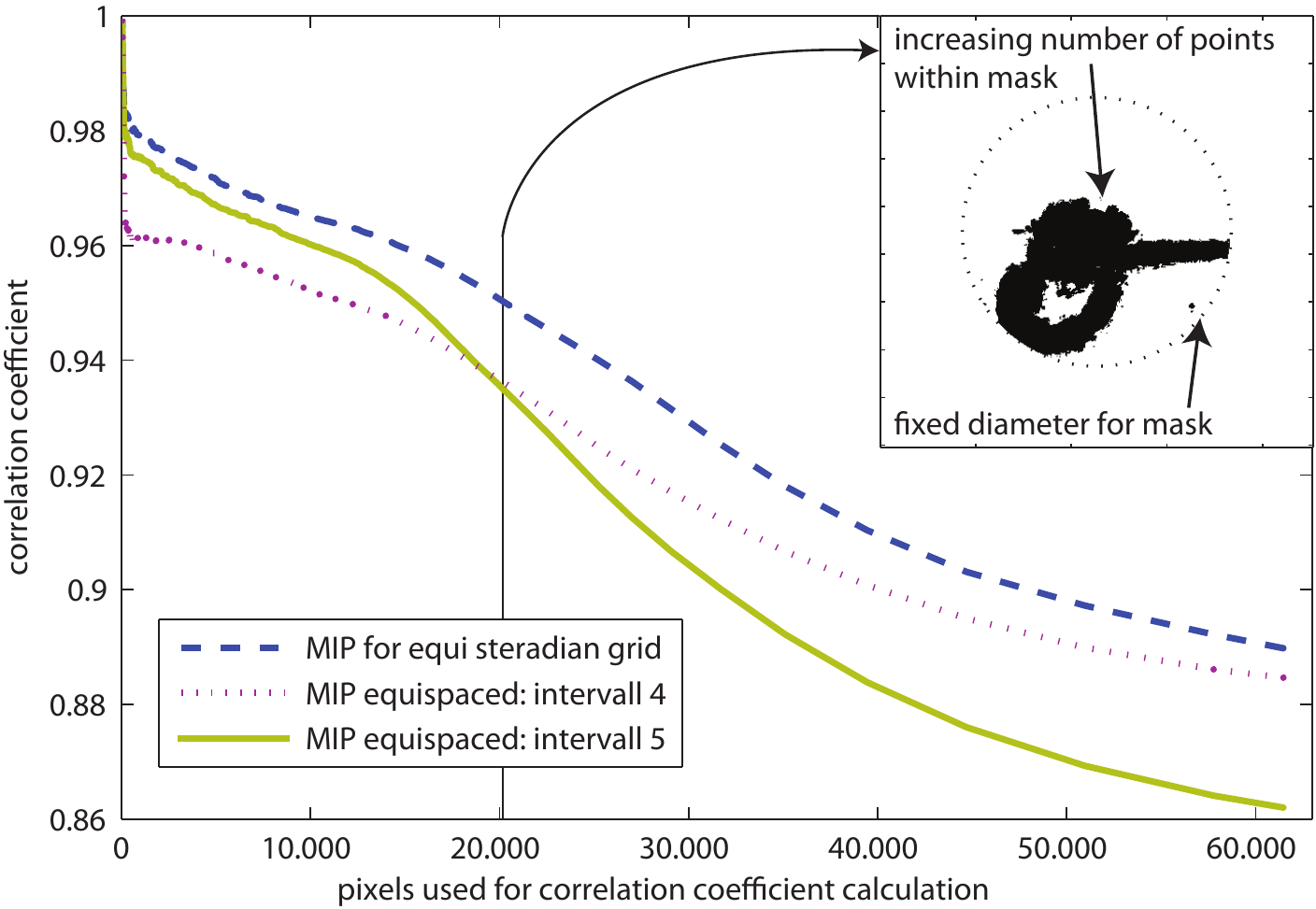}\end{center}

\caption{Correlation coefficient calculated on different number of pixels,
for 3 different reconstructed MIPs. Pixels of interest are added according
to the intensity of the corresponding pixels of the model standard
image (Fig. \ref{fig:Golden-Standard}). In the inlay the correlation
coefficient mask is shown for 20087 points which corresponds to all
pixels with at least 4\% of the maximum value of the model standard and is the value used in Fig. \ref{fig:Correlation-Coefficient-Maskdiameter}. \label{fig:Correlation_Coefficient_mask}}
\end{figure}

In Fig. \ref{fig:Correlation_Coefficient_mask} the equi-steradian
grid is compared with 2 equispaced grids. To get a more precise measure
of the correlation between the reconstructed images and the model
standard, we calculated the correlation coefficient only within a
region of interest. The region of interest is firstly confined by
a centered disc, whose boundary is shown as a dotted circle in the
inlay in Fig. \ref{fig:Correlation_Coefficient_mask}. The $x$-axis
shows the number of pixels of interest used to calculate the correlation
coefficient within this disc. These pixels are increased according
to the intensity of the corresponding pixels of the model's standard
MIP. Fig. \ref{fig:Correlation_Coefficient_mask} demonstrates that
the correlation coefficient for the equi-steradian arrangement, always
remains better, within the depicted disc, when competing against the
two strongest equispaced grids.

\begin{figure}
\begin{center}\includegraphics[width=0.6\columnwidth]{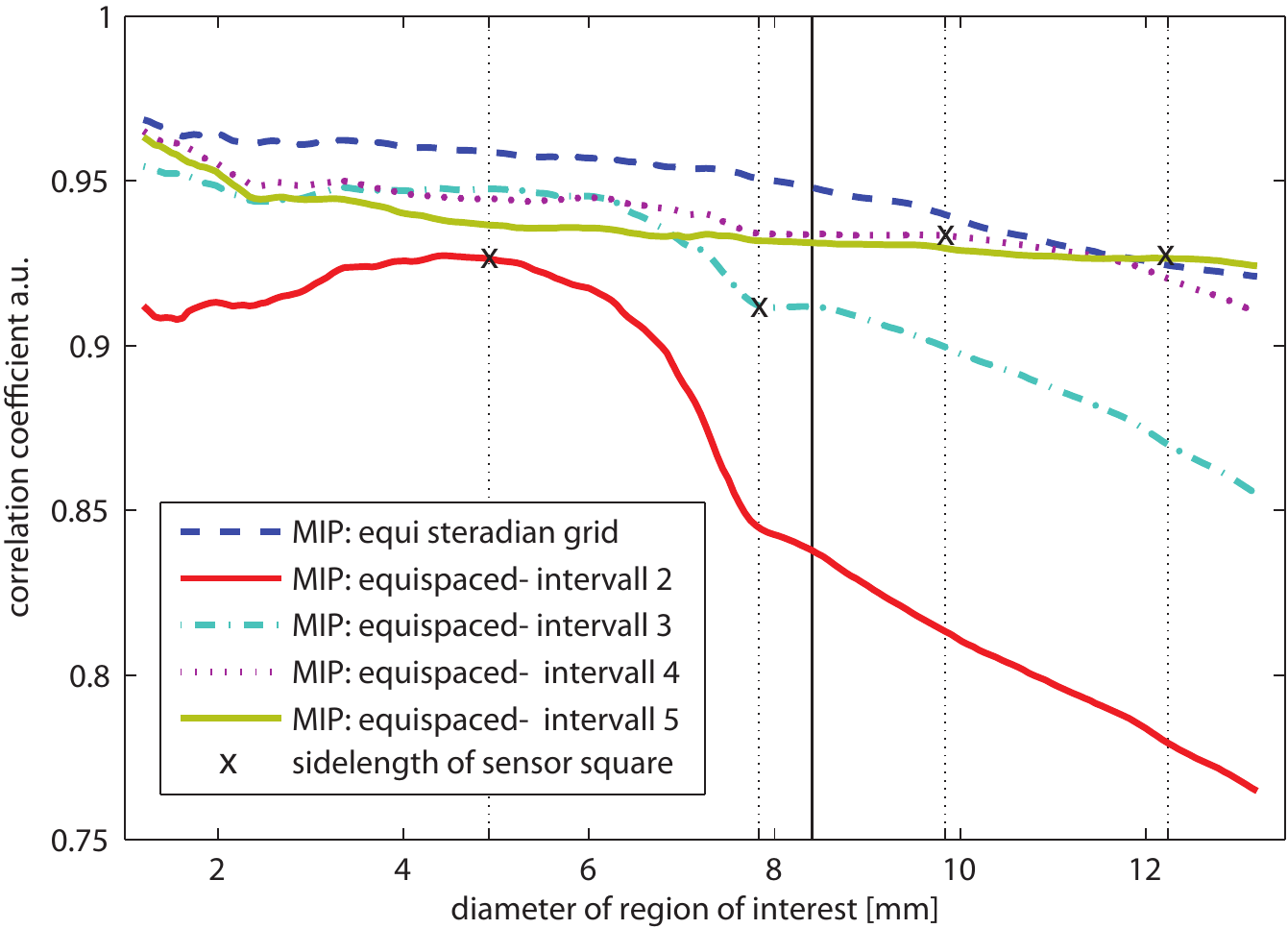}\end{center} 

\caption{Correlation Coefficient calculated on pixels of interest with at least
4\% of the maximum value of the model standard, for an increasing diameter of the confining disc. The straight
line at $8.4\,\mathrm{mm}$ corresponds to the straight line in Fig.
\ref{fig:Correlation_Coefficient_mask}. The other straight lines
indicate the side length of the square of a particular equispaced
sensor point arrangement. \label{fig:Correlation-Coefficient-Maskdiameter}}
\end{figure}

In Fig. \ref{fig:Correlation-Coefficient-Maskdiameter} the 4 equispaced
arrangements, with intervals reaching from 2 to 5 are compared to
the equi-steradian grid. Here the pixels of interest are set to a
threshold of at least 4\% of the maximum value, while the diameter
is increased. The dotted straight lines indicate the side length of
the square of a particular equispaced sensor point arrangement. As
expected the correlation coefficient for this equispaced sensor arrangement
starts to decline around that threshold. 

The correlation coefficient for the equi-steradian grid starts to
fall behind the equispaced grid of interval 5 towards the end. This
is expected, since the equi-steradian grid is only meant to give
better results for a region of interest around the center. This is very clearly the case.
Between a diameter of 2 to 8 mm the correlation coefficient is on
average 25.14 \% closer to a value of $1$ than its strongest equispaced
contender of interval 4, and 30.87 \% better than the interval 5
grid. At a diameter of $7.2\,\mathrm{mm}$ the correlation coefficient
for the equi-steradian grid was 0.960 compared to 0.947 for the interval
4 grid and 0.944 for the interval 5 grid. This is 24.8 \% and 28.6\%
closer to full correlation.

\begin{figure}
\begin{center}\includegraphics[width=0.6\columnwidth]{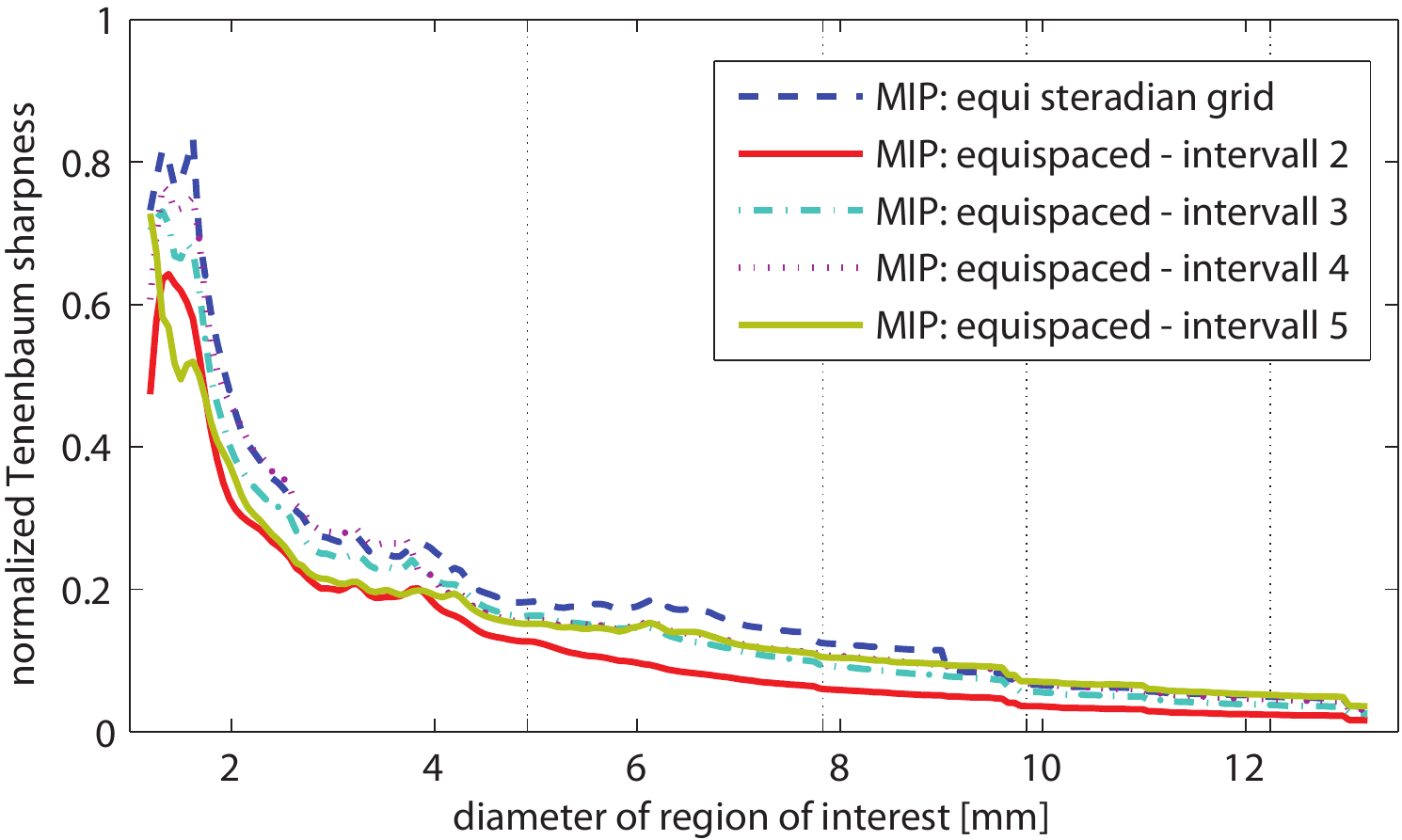}\end{center}\caption{Tenenbaum Sharpness calculated on pixels of interest with at least
4\% of the maximum value of the model standard. The Tenenbaum sharpness
is calculated on the smallest rectangle, that contains all pixels
of interest. \label{fig:Tenenbaum-Sharpness}}
\end{figure}

Fig. \ref{fig:Tenenbaum-Sharpness} shows the normalized Tenenbaum
sharpness. The Tenenbaum sharpness, unlike the correlation coefficient,
cannot be calculated on non-adjacent grid points, therefore it has
been calculated on the smallest rectangle, that contains all pixels
of interest. The equi-steradian has the highest sharpness for most
values, with the interval 4 grid being very slightly better around
a diameter of 3-4 mm. There is a drop of the Tenenbaum sharpness towards
the end.

\section{Discussion and Conclusion}

We computationally implemented a 3D non-uniform FFT
photoacoustic image reconstruction, called NER-NUFFT (non equispaced range-non
uniform FFT) to efficiently deal with the non-equispaced Fourier transform evaluations
arising in the reconstruction formula. This method was compared with the \emph{k-wave} implemented
FFT reconstruction, which uses a polynomial interpolation. 
The two reconstruction methods where compared using 2D targets. The lateral resolution showed an improvement of $18.63\pm 8.5  \%$ 
which is in good agreement with the illustrative results for the star target.
The axial resolution showed an improvement of $168.47 \pm 6.88  \%$.
The computation time was about $30 \%$ less, for the NER-NUFFT than the linearly interpolated FFT reconstruction.
In conclusion the NER-NUFFT reconstruction proved to be unequivocally
superior to conventional linear interpolation FFT reconstruction methods.

We further implemented the NED-NER-NUFFT (non equispaced data-NER-NUFFT),
which allowed us to efficiently reconstruct from data recorded at 
non-equispaced placed sensor points. This newly gained flexibility was used to tackle the
limited view problem, by placing sensors more sparsely further away
from the center of interest. We developed an equiangular sensor placement
for 2D and an equi-steradian placement in 3D,
which assigns one sensor point to each angle/steradian for a given
center of interest. In the 2D computational simulations we showed that this 
arrangement significantly enhances image quality in comparison to regular grids.

In 3D we conducted experiments, where a yarn phantom was recorded.
The maximum intensity projection (MIP) of the full reconstruction
was compared to MIPs of reconstructions that only used about 3\% of
the original data. Within our region of interest, the correlation
of our image was 0.96, which is 24.8\% closer to full correlation
than the best equispaced arrangement, reconstructed from slightly
more sensor points. 

The sensor placement to tackle the limited view problem, combined
with the NED-NER-NUFFT gives significantly better results for an object
located at the center of a bigger sensor surface. This result was
confirmed in the 2D simulation as well as for real data in 3D.

\section*{Acknowledgment}
This work is supported by the Medical University of Vienna, the 
European projects FAMOS (FP7 ICT 317744) and FUN OCT (FP7 HEALTH 201880), 
Macular Vision Research Foundation (MVRF, USA), Austrian Science Fund (FWF), Project P26687-N25 
(Interdisciplinary Coupled Physics Imaging), 
and the Christian Doppler Society 
(Christian Doppler Laboratory "Laser development and their application in medicine").

\def\cprime{$'$}
  \providecommand{\noopsort}[1]{}\def\ocirc#1{\ifmmode\setbox0=\hbox{$#1$}\dimen0=\ht0
  \advance\dimen0 by1pt\rlap{\hbox to\wd0{\hss\raise\dimen0
  \hbox{\hskip.2em$\scriptscriptstyle\circ$}\hss}}#1\else {\accent"17 #1}\fi}
  \def\cprime{$'$}

\section{Appendix: Algorithm for equi-steradian sensor arrangement}
\label{App:equi-ster}

In our algorithm, the grid size and the distance of the center of
interest from the sensor plane is defined. The number of sensor points
$N$ will be rounded to the next convenient value.

Our point of interest is placed at $z=r_{0}$, centered at a square
$xy$ grid. The point of interest is the center of a spherical coordinate
system, with the polar angle $\theta=0$ at the $z$-axis towards
the $xy$-grid. 

First we determine the steradian $\Omega$ of the spherical cap from
the point of interest, that projects onto the acquistion point plane via 
\[\Omega=2\pi\left(1-\cos\left(\theta_{max}\right)\right)\;.\]
%\commentO{Achtung, das $\Omega$ is was anderes als auf Seite 3, oder?}
This leads to a unit
steradian $\omega=\Omega/N$ with $N$ being the number of sensors
one would like to record the signal with. The sphere cap is then subdivided
into slices $k$ which satisfy the condition 
\[\omega\, j_{k}=2\pi\left(\cos\left(\theta_{k-1}\right)-\cos\left(\theta_{k}\right)\right)\,,\]
where $\theta_1$ encloses exactly one unit steradian $\omega$ and $j_k$ has to be a power of two, in order to guarantee some symmetry. 
The value of $j_k$ doubles, when $r_s>1.8\cdot r_k$, where $r_s$ is the chord length between two points on $k$ and $r_k$ is the distance to the closest point on $k-1$.
These values are chosen in order to guarantee 
We apply some restrictions to approximate equidistance between acquisition points on the sensor surface.

The azimuthal angles for a slice $k$ are calulated according to:
\[\varphi_{i,k}=\left(2\pi i\right)/j_{k}+\pi/j_{k}+\varphi_{r}\,,\]
with $i=0,\ldots,j-1$, where 
\[\varphi_{r}=\varphi_{j_{k-1},k-1}+\left(k-1\right)2\pi/(j_{k-1})\]
stems from the former slice $k-1$ . 
The sensor points are now placed on the $xy$-plane at the position indicated
by the spherical angular coordinates:
\[\left(\mathrm{pol,az}\right)=\left(\left(\theta_{k}+\theta_{k+1}\right)/2,\varphi_{i,k}\right)\]

\section{Appendix: Quality measures}
\label{App:Qualitmeas}

The correlation coefficient $\rho$ is a 
measure of the linear dependence between two images $U_{1}$ and $U_{2}$.

Its range is $\left[-1,1\right]$ and a correlation coefficient close
to 1 indicates linear dependence \cite{Sch11}. It is defined via
the variance, $\mbox{Var}\left(U_{i}\right)$ of each image and the
covariance, $\mbox{Cov}\left(U_{1},U_{2}\right)$ of the two images:
\begin{equation}\label{eq:corrcoeff}
\begin{array}{cl}
\rho\left(U_{1},U_{2}\right)=\frac{\mbox{Cov}\left(U_{1},U_{2}\right)}{\sqrt{\mbox{Var}\left(U_{1}\right)\mbox{Var}\left(U_{2}\right)}}\;.
\end{array}
\end{equation}

We did not choose the widely used $L^{p}$ distance measure because
it is a morphological distance measure, meaning it defines the distance
between two images by the distance between their level sets. Therefore
two identical linearly dependent images can have a correlation coefficient
of $1$ and still a huge $L^{p}$ distance. Normalizing the images only mitigates this problem,
because single, high intensity artifacts or a varying intensity over
the whole image can greatly alter the $L^{p}$ distance.

While the correlation coefficient is a good measure for the overall
similarity between two images it does not include any sharpness measure.
Hence blurred edges are punished very little, in comparison to slight variations
of homogeneous areas. To address this shortcoming we chose a sharpness
measure or focus function as a second quality criterion. Sharpness
measures are obtained from some measure of the high frequency content
of an image \cite{GroYouLig85}. They have also been used to select
the best sound speed in photoacoustic image reconstruction \cite{TreVarZhaLauBea11}.
Out of the plethora of published focus functions we select the Tenenbaum
function, because of its robustness to noise:

\begin{equation}
\begin{array}{cl}
F_\text{Tenenbaum}=\underset{x,y}{\sum}\left(g*U_{x,y}\right)^{2}+\left(g^{T}*U_{x,y}\right)^{2}\,,
\end{array}
\end{equation}

with $g$ as the Sobel operator: 
\begin{equation}
g=\left(\begin{array}{ccc}
-1 & 0 & 1\\
-2 & 0 & 2\\
-1 & 0 & 1
\end{array}\right)\;.
\end{equation}

Like the $L^{2}$ norm and unlike the correlation coefficient the
Tenenbaum function is an extensive measure, meaning it increases with
image dimensions. Therefore we normalized it to $\overline{F}_\text{Tenenbaum}=F_\text{Tenenbaum}/N$,
where $N$ is the number of elements in $U$. 

\end{document}